%% file: buscount11-Jun_08_macs_inc.tex
\documentclass[10pt]{article}
\usepackage{amssymb}
\usepackage{amsmath}
\usepackage{amsthm}
\usepackage{epsfig}

\newcommand{\gw}{G_{\scriptscriptstyle{W}}}
\newcommand{\phit}{\varphi_{\tau}}

\newcommand{\phiw}{\varphi_{\scriptscriptstyle{W}}}

\newcommand{\be}{\begin{eqnarray}}
\newcommand{\ee}{\end{eqnarray}}

\newcommand{\R}{{\bf R}}

\newcommand{\dist}{{\rm dist}\,}

\newcommand{\cof}{{\rm cof}\,}

\newcommand{\tbyt}{\mathbf{R}^{2 \times 2}}

\newcommand{\sct}{\mathcal{T}}
\newcommand{\scs}{\mathcal{S}}

\newcommand{\romx}{\rho_{\scriptscriptstyle{\textrm{max}}}}
\newcommand{\romn}{\rho_{\scriptscriptstyle{\textrm{min}}}}
\newcommand{\scl}{\mathcal{L}}

\newtheorem{thm}{Theorem}[section]
\newtheorem{prop}[thm]{Proposition}
\newtheorem{lem}[thm]{Lemma}

\newtheorem{rem}[thm]{Remark}
\numberwithin{equation}{section}
\title{A remark on the structure of the Busemann representative of a polyconvex function}
\author{J. J. Bevan
\footnote{Address for correspondence: Department of Mathematics, University of Surrey, Guildford, GU2 7XH, UK.}}
\date{\today}
\begin{document}
\maketitle

\section{Introduction}\label{intro}

Polyconvexity was first identified by Morrey in \cite{Mo52} and was
later developed by Ball \cite{jb1} in connection with nonlinear
elasticity.   A function $W:\R^{N \times n} \to \R \cup \{\infty\}$
is polyconvex if there exists a convex function $\varphi$, said to be a
convex representative of $W$, such that
\[W(\xi)=\varphi(R(\xi))\]
for all real $N \times n$ matrices $\xi$, where $R(\xi)$ is the list of minors of $\xi$ written in some fixed order.  Busemann \emph{et al.} pointed out in \cite{Bu1} that there is a largest such convex representative: we refer to this as the Busemann representative and denote it by $\phiw$.

One of the broader aims of the series of papers \cite{Bu1} Busemann \emph{et al.} was to study the restriction of convex functions to non-convex sets.  Ball observed in \cite{jb1} that polyconvexity fits into this framework, and the relationship between the two has since been explored further in \cite{busrep}.

The Busemann representative $\phiw$ of a given polyconvex function $W:\R^{N \times n} \to \R \cup \{\infty\}$ can be expressed as 
\begin{equation}\label{lagrange}\phiw(X)=\inf\left\{\sum_{i=1}^{d+1}\lambda_{j}W(\xi_{j}): \ \ \lambda_{j} \geq 0, \sum_{j=1}^{d+1}\lambda_{j}=1 \ \textrm{and} \ \sum_{j=1}^{d+1}\lambda_{j}R(\xi_{j})=X\right\}.
\end{equation}
\noindent Here, $d$ is the least integer such that $R(\xi) \in \R^{d}$ for all $\xi \in \R^{N \times n}$ and $X$ lies in $\R^{d}$.  Busemann \emph{et al.}
proved that 
\begin{equation}\label{null}\phiw(X)= \sup \{a(X): a \in \scl\},\end{equation} 
where 
\[\scl = \{\phi \ \textrm{affine}: \ \phi(R(\xi)) \leq W(\xi) \ \forall \ \xi \in \R^{N \times n} \}.\]
The graph of any $\phi \in \scl$ is a hyperplane, so \eqref{null} states that $\phiw$ is built from hyperplanes which lie below the set 
$G_{\scriptscriptstyle{W}}:=\{(R(\xi), W(\xi)): \ \xi \in \R^{N \times n}\}$.

The main result in this short note is that there is no redundancy in the expression \eqref{null} in the case $N=n=2$.  To be precise, one cannot replace $\scl$ in \eqref{null} by the smaller class $\sct$, where
\[\sct=\{\phi \in \scl: \exists \ \xi \in \tbyt \ \textrm{s.t.} \ W(\xi)=\phi(R(\xi))\}.\]  Thus $\sct$ represents the collection of supporting hyperplanes to $\gw$ which meet $\gw$ in at least one point.   We define
\[\phit(X)= \sup \{a(X): \ a \in \sct\}.\]
Note that $\phiw \geq \phit$ in view of the inclusion $\sct \subset \scl$.
It is proved in the next section that for a certain choice of $W$ it is the case that $\phiw > \phit$ on a large set.  
This result is surprising since the set $\{R(\xi): \ \xi \in \tbyt\}$ is large: its convex hull is the whole of $\R^{5}$. (For a proof of this fact see \cite{jb1}.)  Certainly one might expect $\phiw=\phit$ to be the case under extra assumptions, which could include super-quadratic growth of $W$, for example.  See \cite{busrep} for further details. 

The result of this note is relevant to \cite[Lemma 2.4]{busrep}, where the structure of $\phiw$ is described.   We present a version of the lemma below for the reader's benefit; for the proof consult \cite{busrep}.

\begin{lem}{\cite[Lemma 2.4]{busrep}}
Let $\scs=\{R(\xi): \ \xi \in \R^{N \times n}\}$ and suppose $W:\R^{N \times n} \to \R$ is polyconvex.
Define $\phiw$ by \eqref{lagrange}.
Then for each $X \in \R^{d}$ either one or both of the following hold:
\begin{itemize}\item[(a)]there exists $Y \in \scs$ such that
    $\phiw\arrowvert_{[Y,X]}$ is affine;
\item[(b)]there exists a unit vector $e \in \R^{d}$ such that for all
  $Y \in \R^{d}$ and all $t \in \R$ the function $t \mapsto
  \phiw(Y+te)$ is constant.
\end{itemize}
\end{lem}

The dichotomy can be sharp in the sense that (a) and not (b) can hold, as easy examples show, and that (b) and not (a) can hold, which is a consequence of the counterexample constructed below.  It is shown in \cite{busrep} that when (a) holds the differentiability of $\phiw$ on $\scs$ implies that $\phiw$ is the unique convex representative.  The counterexample below shows that this result is false when (b) holds and (a) does not.

\subsection{Notation}  We do not distinguish between the inner product of two matrices and the inner product of two vectors in $\R^{5}$, using $\cdot$ for both.  Here,   $\R^{5}$ is shorthand for $\tbyt \times \R$, and in this case the inner product of $(\xi, s)$ with $(\eta, t)$ is given by
\[(\xi, s) \cdot (\eta, t) = \xi \cdot \eta + st,\]
where $\xi, \eta$ are two matrices in $\R^{2 \times 2}$, $s,t \in \R$ and  
\[\xi \cdot \eta = \textrm{tr}(\xi^{T} \eta).\]
Finally, if $a,b \in \R^{2}$ then the $2 \times 2$ matrix $a \otimes b$ has $(i,j)-$entry $a_{i}b_{j}$.  

\section{Construction of $W$ such that $\phiw  > \phit$ on a large set}

We restrict attention to polyconvex functions defined on $\R^{2 \times 2}$, so that $R(\xi)=(\xi,\det \xi)$ for each $2 \times 2$ matrix $\xi$.  To begin with we recall some basic facts about the subgradients of $\phiw$ (for the definition of the subgradient of a convex function see \cite{Rock}).  When $W:\R^{2 \times 2} \to \R$ is polyconvex and differentiable on an open set $U \subset \R^{2 \times 2}$ it can be shown that for each $\xi \in U$
\begin{equation}\label{part}\partial \phiw (R(\xi))=\{(DW(\xi)-\rho \cof \xi, \rho): \  \romn(\xi) \leq \rho \leq \romx(\xi)\},\end{equation}
where the functions $\romx$, $\romn: \R^{2 \times 2} \to \R$ are defined 
by 
\begin{eqnarray}\label{rows}\romx(\xi)&=&\inf \left\{\frac{W(\eta+\xi)-W(\xi)-DW(\xi) \cdot \eta}
{\det \eta}: \  \det \eta > 0\right \} \\
\romn(\xi)&=&\sup \left\{\frac{W(\eta+\xi)-W(\xi)-DW(\xi) \cdot \eta}
{\det \eta}: \  \det \eta< 0\right \}.
\end{eqnarray}
The proof of these assertions can be found in \cite[Section 2]{busrep}.  Thus when $\xi \in U$, a sufficient condition for the differentiability of $\phiw$, and hence of $\phit$ (because $\phiw \geq \phit$ on $\R^{5}$, and because $\phiw$ and $\phit$ agree on $\scs$---see \cite[Corollary 2.5]{bkk}), at $R(\xi)$ is that there exists a number $\rho(\xi)$ such that 
\[W(\xi+\eta) \geq W(\xi)+DW(\xi)\cdot \eta + \rho(\xi) \det \eta\]
for all $2 \times 2$ matrices $\eta$.  

Now let $[\xi]=\xi-\xi_{11}e_{1} \otimes e_{1}$, where $e_{1}$ is the first canonical basis vector in $\R^{2}$, and define $W(\xi)=|([\xi], \det \xi - y)|$, where $|z|$ is the usual Euclidean norm in $\R^{5}$ and where $y$ is a fixed positive number.   It is easy to see that $W$ is polyconvex and differentiable away from the set 
$\{\xi: \ W(\xi)=0\}$, which, since $y \neq 0$, is empty.    
With the above remarks in mind the following proposition shows that $\phiw$ is differentiable at all points $\R(\xi)$ in 
$\scs$.

\begin{prop}\label{half} Let $\xi \in \tbyt$ and let $W$ be as above.  Then for all $\eta$
\[ W(\xi+\eta)-W(\xi)-DW(\xi)\cdot \eta \geq \rho(\xi)\det \eta , \]
where $\rho(\xi)=\frac{(\det \xi -y)}{W(\xi)}$.
\end{prop}
\begin{proof}
The inequality amounts to proving
\[|([\xi + \eta],\det(\xi+\eta)-y)| \geq \frac{1}{W(\xi)}([\xi + \eta] \cdot [\xi] +(\det \xi -y)(\det (\xi + \eta)-y)).\]
  But this follows directly from the Cauchy-Schwarz inequality.
\end{proof}
 
\begin{rem}\emph{The choice of $\rho(\xi)$ in Proposition \ref{half} is by analogy with the following example.  Suppose $f(\xi)=|R(\xi)|$ and note that an obvious convex representative of $f$ is $\varphi(\xi,\delta)=|(\xi,\delta)|$.  Differentiating this with respect to $\delta$, evaluating at $R(\xi)$, where $\xi \neq 0$, and referring to \eqref{part} gives a candidate $\rho(\xi)=\frac{\det \xi}{f(\xi)}$. }
\end{rem}

\begin{figure}[ht]
\centering
\input{tuppence.pstex_t} 
\caption{A graphical interpretation of the contructions of $W$, $\phiw$ and $\phit$.  $\scs$ can be thought of as the union of the two curves in the $x-y$ plane, the graph of $\phiw$ as the union of the plane $ABCD$ together with the two sloping planes it meets at $AD$ and $BC$, and the graph of $\phit$ as the union of the two sloping planes.  The function $W$ is represented by the restriction of $\phiw$ to $\scs$; its graph is shown with dotted lines.  Clearly, $\phiw > \phit$ in the projection of $ABCD$ in the $x-y$ plane.}
\label{tupp}
\end{figure}
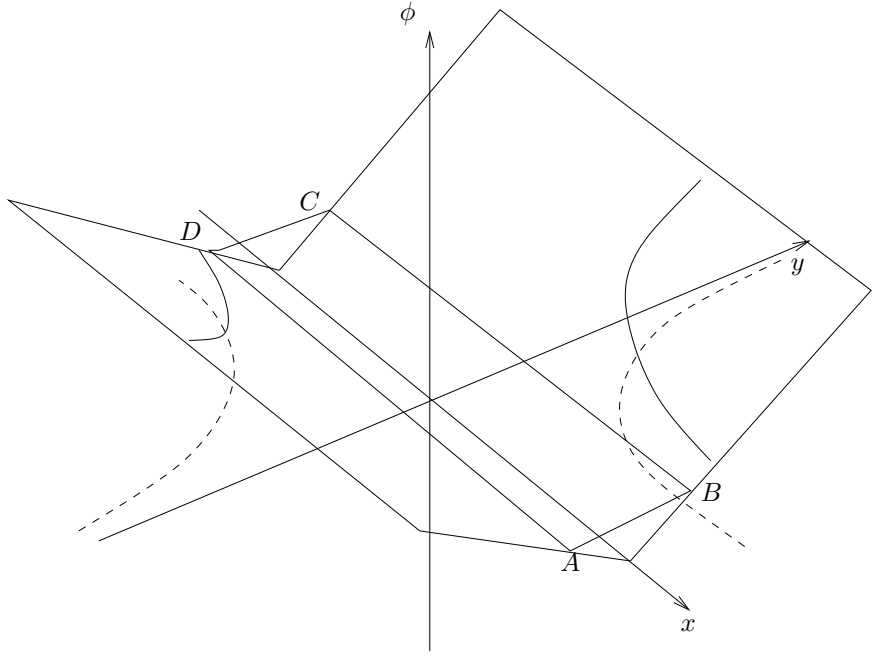

Now consider the line $L:=\textrm{Span}\{e_{1} \otimes e_{1}\}$.  Clearly $\det l = 0$ for all $l \in L$.  Since 
$D^{2}\det( \xi)[\eta,\eta]=2\det \eta$ for all $2 \times 2$ matrices $\xi$ and $\eta$, we can assume that the curvature of the
graph of the determinant (i.e., the curvature of $\scs$) is bounded above uniformly on the set $\{l+\eta: l \in L, \ |\eta| < 1\}$.  In particular, we deduce that for sufficiently small $\epsilon > 0$ the (convex) tube 
\[T_{\epsilon}:=\{(l+\eta,y): \ l \in L, \ |\eta|\leq \epsilon\},\] 
which lies in $\R^{5}$, satisfies $\dist(T_{\epsilon},\scs)>0$.   With $W$ as above it is claimed that $\phiw > \phit$ on the tube $T_{\epsilon}$.  Figure \ref{tupp} below is intended as an analogy which may help the reader to visualize the idea behind the proof of Proposition \ref{fan}.

\begin{prop}\label{fan}Let $W(\xi)=|([\xi],\det \xi - y)|$ and assume $\epsilon$ has been chosen so that the tube $T_{\epsilon}$ does not meet $\scs$.  Then $\phiw(X)>\phit(X)$ for all $X \in T_{\epsilon}$.
\end{prop}
\begin{proof}  Recall that $\phit(X)=\sup\{a(X): a \in \sct\}$, where $\sct$ consists of all those affine functions $a$ satisfying $a(\xi,\det\xi) \leq W(\xi)$ for all $\xi \in \tbyt$, and $a(\xi_{0},\det \xi_{0})=W(\xi_{0})$ for at least one $\xi_{0}$.  Suppose $a_{\xi_{0}}$ is such that $a_{\xi_{0}}(\xi_{0},\det \xi_{0})=W(\xi_{0})$.
Standard arguments from convex analysis together with the differentiability of $\phiw$ (Proposition \ref{half} above) at all $(\xi_{0},\det \xi_{0})$ show that the gradient of the affine function $a_{\xi_{0}}$ at $(\xi_{0},\det \xi_{0})$ must be 
$D\phiw(\xi_{0},\det \xi_{0})$.  Since $a_{\xi_{0}}$ is affine, and in view of \eqref{part},  it follows that for all $X$ in $\R^{5}$
\begin{eqnarray*}a_{\xi_{0}}(X) & = & W(\xi_{0})+D \phiw (\xi_{0},\det \xi_{0}) \cdot (X-(\xi_{0},\det \xi_{0})) \\
& = & ([\hat{X}],X'-y)) \cdot \frac{([\xi_{0}], \det \xi_{0} - y)}{W(\xi_{0})}.
\end{eqnarray*}
Here we have used the notation $X=(\hat{X},X') \in \tbyt \times \R$.  Thus 
\begin{equation}\label{final}\phit(X)=\sup \left\{([\hat{X}],X'-y)) \cdot \frac{([\xi_{0}], \det \xi_{0} - y)}{W(\xi_{0})}: \ \xi_{0} \in \tbyt\right\}.\end{equation}
Provided we can find $\xi_{0}$ such that $([\hat{X}],X'-y))$ and $([\xi_{0}], \det \xi_{0} - y)$ are parallel, or asymptotically parallel (which will be made clear below), then it will follow essentially from the Cauchy-Schwarz inequality that $\phit(X)=|([\hat{X}],X'-y))|$.
There are three cases to consider, and in doing so we shall refer to the unit vector 
$\frac{([\xi_{0}], \det \xi_{0} - y)}{W(\xi_{0})}$ by $u(\xi_{0})$.

\vspace{1mm}
\noindent\textbf{(i)} $[\hat{X}]=0$.  Note that $u(0)=(0,-1)$, which gives $\phit(X)=|X'-y|$ provided $y > X'$.  Otherwise note that $u(k Q) \to (0,1)$ as $k \to \infty$ whenever $Q$ is a rotation matrix (i.e. $Q \in SO(2)$), which implies $u(kQ) \cdot (0,X'-y) \to |X'-y|$ whenever $X' > y$.  If $y = X'$ then $\phit(x') = |X' - y|$ trivially.

\vspace{1mm}
\noindent\textbf{(ii)} $\hat{X}_{22}\neq 0$.  Set $\xi_{0}=[\hat{X}]$ and consider $\xi_{\mu}=\xi_{0}+\mu e_{1} \otimes e_{1}$.  We require $\det \xi_{\mu}=X'$.  But this can easily be satisfied by an appropriate choice of $\mu$, and on using $\hat{X}_{22} \neq 0$ in $\det \xi_{\mu} = \det \xi_{0} + \mu \hat{X}_{22}$.

\vspace{1mm}
\noindent\textbf{(iii)} $[\hat{X}]\neq 0,  \ \hat{X}_{22}=0$.   As before, choose $\xi_{0}$ to satisfy $\xi_{0}=[\hat{X}]$ and let $\xi_{\mu,\nu}=\xi_{0}+\mu e_{1} \otimes e_{1} + \nu e_{2} \otimes e_{2}$, where $\mu$ and $\nu$ are parameters. Now we seek $\mu$ and $\nu$ such that $\det \xi_{\mu,\nu}=X'$, that is,  
\begin{equation}\label{peacock}\mu \nu=X'+\hat{X}_{12}\hat{X}_{21}.\end{equation}
But  $[\xi_{\mu,\nu}]=[\hat{X}]+\nu e_{2} \otimes e_{2}$, and hence 
\[u(\xi_{\mu,\nu}) \to \frac{([\hat{X}],X'-y)}{|([\hat{X}],X'-y)|}\]
provided $\mu \to \infty$ and $\nu \to 0$ consistent with \eqref{peacock}.  

Thus in each case we have $\phit(X)=|([\hat{X}],X'-y)|$.  To conclude the proof note that $W(\xi)$ can be interpreted as the distance of the point $(\xi,\det \xi)$ to the centre of the tube $T_{\epsilon}$.  The construction of $T_{\epsilon}$ above therefore implies that 
$W(\xi) \geq \epsilon$ for all $2 \times 2$ matrices $\xi$.  Hence $\phiw(X) \geq \epsilon$ for all $X$, while $\phit(X) < \epsilon$ for all $X$ inside the tube $T_{\epsilon}$.
 
\end{proof}

With reference to the statement of \cite[Lemma 2.4]{busrep} given in the introduction, we remark that because alternative (a) of \cite[Lemma 2.4]{busrep} fails for points $X$ in the tube $T_{\epsilon}$ it must be that (b) holds for such points.   It was shown in  \cite[Proposition 3.5]{busrep} that if alternative (a) held at all $X$ and if $\phiw$ was differentiable on $\scs$ then $\phiw$ was the unique convex representative  of $W$.   This result is clearly false when alternative (b) holds at some $X$, even when, as we have seen above, $\phiw$ is differentiable on $\scs$.

\vspace{3mm}
\noindent \textbf{Acknowledgement}  This research was supported by an EPSRC Postdoctoral Research Fellowship  GR/S29621/01, by the European Research and Training Network MULTIMAT and by an RCUK Academic Fellowship.  I thank Prof. B. Kirchheim for reading a draft version of the paper and for the idea leading to Figure \ref{tupp}.

\end{document}

%% file: tuppence.pstex_t
\begin{picture}(0,0)%
\includegraphics{tuppence.pstex}%
\end{picture}%
\setlength{\unitlength}{3315sp}%
\begingroup\makeatletter\ifx\SetFigFont\undefined%
\gdef\SetFigFont#1#2#3#4#5{%
  \reset@font\fontsize{#1}{#2pt}%
  \fontfamily{#3}\fontseries{#4}\fontshape{#5}%
  \selectfont}%
\fi\endgroup%
\begin{picture}(6474,4893)(2314,-5923)
\put(6451,-5311){\makebox(0,0)[lb]{\smash{{\SetFigFont{10}{12.0}{\rmdefault}{\mddefault}{\updefault}$A$}}}}
\put(7501,-4786){\makebox(0,0)[lb]{\smash{{\SetFigFont{10}{12.0}{\rmdefault}{\mddefault}{\updefault}$B$}}}}
\put(4501,-2611){\makebox(0,0)[lb]{\smash{{\SetFigFont{10}{12.0}{\rmdefault}{\mddefault}{\updefault}$C$}}}}
\put(3601,-2836){\makebox(0,0)[lb]{\smash{{\SetFigFont{10}{12.0}{\rmdefault}{\mddefault}{\updefault}$D$}}}}
\put(7351,-5761){\makebox(0,0)[lb]{\smash{{\SetFigFont{10}{12.0}{\rmdefault}{\mddefault}{\updefault}$x$}}}}
\put(8176,-3061){\makebox(0,0)[lb]{\smash{{\SetFigFont{10}{12.0}{\rmdefault}{\mddefault}{\updefault}$y$}}}}
\put(5251,-1186){\makebox(0,0)[lb]{\smash{{\SetFigFont{10}{12.0}{\rmdefault}{\mddefault}{\updefault}$\phi$}}}}
\end{picture}%